\documentclass[leqno,12pt]{amsart} 
\usepackage{amssymb} 
\theoremstyle{plain} 
\newtheorem{theorem}{\indent\sc Theorem}[section]
\newtheorem{lemma}[theorem]{\indent\sc Lemma}

\theoremstyle{definition} 

\usepackage{url}
\usepackage{amsthm}
\usepackage{amssymb}
\usepackage{latexsym}
\usepackage{amsfonts}
\usepackage{amsmath}
\usepackage{amstext}
\usepackage{accents}
\usepackage{cancel}
\usepackage{eucal}
\usepackage{amscd}
\usepackage{hyperref}
\usepackage[shortalphabetic]{amsrefs}
\usepackage{color}
\usepackage{pgfpages}
\usepackage{graphicx}
\usepackage[utf8]{inputenc}
\usepackage[T1]{fontenc}
\usepackage{animate}
\usepackage{tikz}

\newcommand{\R}{\mathbb{R}}

\renewcommand{\vec}[1]{\ensuremath{{\bf #1}}}

\begin{document}

\title[Ideal surfaces with boundary]{A rigidity theorem for ideal surfaces with flat boundary}

\author[J. McCoy]{James McCoy}
\address{School of Mathemtical and Physical Sciences, University of Newcastle}
\email{James.McCoy@newcastle.edu.au}

\author[G. Wheeler]{Glen Wheeler}
\address{Institute for Mathematics and its Applications, University of Wollongong}
\email{glenw@uow.edu.au}


\thanks{Research supported by the Australian Research Council, Discovery Project DP150100375. The authors also acknowledge Benjamin Maldon (University of Newcastle) for assistance with typesetting through the University of Newcastle Priority Research Centre for Computer-Assisted Research Mathematics and its Applications (CARMA)}

\begin{abstract}
We consider surfaces with boundary satisfying a sixth order nonlinear elliptic partial differential equation corresponding to extremising the $L^2$-norm of the gradient of the mean curvature.  We show that  such surfaces with small $L^2$-norm of the second fundamental form and satisfying so-called `flat boundary conditions' are necessarily planar.
\end{abstract}

\keywords{higher order geometric partial differential equation, sixth order elliptic equation, Neumann boundary condition}
\subjclass[2010]{35J30, 58J05, 35J62}
\maketitle

\section{Introduction} \label{S:intro}
\newtheorem{main}{Theorem}[section]
We are interested extremal surface of the geometric energy
\begin{equation} \label{E:energy}
  \mathcal{F}[ f] =\int_\Sigma \left| \nabla H \right|^2 d\mu
  \end{equation}
under the hypothesis
\begin{equation} \label{E:smallness}
\int_\Sigma \left| A \right|^2 d\mu \leq \varepsilon_0
\end{equation}
where $\varepsilon_0>0$ is a small, universal constant.  Here $f: \Sigma \rightarrow \mathbb{R}^3$ a smooth immersion of surface $\Sigma$ with boundary; $d\mu$ is the induced surface area element; $H = \kappa_1 + \kappa_2$ and $\left| A\right|^2= \kappa_1^2 + \kappa_2^2$ are respectively the mean curvature and the norm of the second fundamental form of $f\left( \Sigma \right)$ and $\nabla$ is the covariant derivative on $f\left( \Sigma \right)$.  Such extremal surfaces we will call `ideal surfaces'.  Clearly minimal surfaces and surfaces of constant mean curvature are absolute minimisers of \eqref{E:energy}.

Previous work related to the type of result of this article includes rigidity of manifolds with Ricci curvature bounds whose volume is close to that of the sphere \cite{A90} and work on rigidity and classification of minimal submanifolds \cite{CDCK, FC, L69, S68, R70}, and hypersurfaces of constant mean curvature \cite{NS69}.  More recently and closer in spirit to our result here is work on higher-order geometric gap lemmas.  The first of these, for Willmore surfaces, appeared in \cite{KS} (Theorem 2.7) as part of a larger investigation of the gradient flow for the Willmore functional
$$\int_\Sigma \left| \textbf{H} \right|^2 d\mu$$
for surfaces $\Sigma$ without boundary immersed in $\mathbb{R}^{n+k}$, $k\geq 1$, satisfying the smallness condition 
\begin{equation} \label{E:smallA0}
  \int_\Sigma \left| A^0 \right|^2 d\mu \leq \varepsilon_0 \mbox{.}
\end{equation}
Here $\left| A^0 \right|^2 = \left| A \right|^2 - \frac{1}{2} H^2 = \frac{1}{2} \left( \kappa_1 - \kappa_2\right)^2$, the norm of the trace-free second fundamental form, is a pointwise measure of closeness to each other of the two principal curvatures $\kappa_1$ and $\kappa_2$.  The gap lemma of \cite{KS} gives that the resulting surface without boundary is either an embedded plane or sphere.  Later the second author of this article obtained a gap lemma for stationary solutions of the surface diffusion flow \cite{W12} without boundary and under the smallness condition
 $$\int_\Sigma \left| A^0 \right|^n d\mu \leq \varepsilon_0$$
 for surfaces of dimension $n=2$ and $n=3$; again such a surface is either an embedded plane or sphere.  The second author further obtained a gap lemma for biharmonic surfaces in \cite{W13} and, together with the first author, for some Helfrich surfaces \cite{MW}.  This result was extended to some other Helfrich surfaces in \cite{BWW}.  With Parkins the authors obtained a gap lemma for triharmonic surfaces \cite{MPW}; Parkins additionally obtained the corresponding result for polyharmonic surfaces in \cite{P17}.
 
 For several of the abovementioned results there are also versions for surfaces with boundary, with either of two boundary conditions:
\begin{enumerate}
	\item \emph{umbilic boundary conditions} $\left| \nabla  A^0 \right| = \left| A^0 \right| = 0$; or
	\item \emph{flat boundary conditions} $\left| \nabla A  \right| = \left| A \right| = 0$.
\end{enumerate}
With suitable smallness conditions, umbilic boundary conditions lead to parts of spheres and planes, while flat boundary conditions allow parts of planes only \cite{W14}.  

In many cases above results hold for arbitrary codimension.\\

The main result of this article may be stated as follows:

\begin{main} \label{T:main}
Suppose $f: \Sigma \rightarrow \mathbb{R}^3$ satisfies 
\begin{equation} \label{E:theeqn}
\mathcal{I}[f]:=  \Delta^2 H + \left| A\right|^2 \Delta H - \left( A^0\right)^{ij} \nabla_i H \nabla_j H = 0 
\end{equation}
with boundary conditions 
\begin{equation} \label{E:BCs}
\left| A \right| = 0 \mbox{ and } \nabla_\eta H = \nabla_\eta \Delta H =0 \mbox{.}
\end{equation}
	If $f$ also satisfies \eqref{E:smallness}
	for $\varepsilon_0>0$ sufficiently small, then $f\left( \Sigma \right)$ is part of a flat plane.
\end{main}

Above we have used $\eta$ to denote the unit conormal to the boundary and $\Delta$ the Laplace-Beltrami operator.  We also have $\Delta^2 H := \Delta \Delta H$ and use the standard Einstein summation convention of summing over repeated indices.\\

\noindent \textbf{Remarks:} 
\begin{enumerate}
  \item[1.] The boundary conditions \eqref{E:BCs} are understood in the sense of limits approaching the boundary within the surface.  We refer the reader to \cite{W14} for precise definitions.
  \item[2.] In the case of $f: \Sigma \rightarrow \mathbb{R}^{k}$, $k>3$, \eqref{E:theeqn} may be replaced by relatively weak orthogonality condition
\begin{equation*} 
  \left< \mathcal{I}\left[ f\right], \vec{H} \right>_{\mathbb{R}^k} = 0 \mbox{.}
\end{equation*}
With analogous boundary conditions \ref{E:BCs} and smallness condition \ref{E:smallness} we obtain the same result as Theorem \ref{T:main}.  In this article we restrict to the case $k=3$ for notational simplicity; the workings in the general case are essentially the same. 
 \item[3.] In the case $k=3$ \eqref{E:theeqn} may be replaced by the slightly more general
 $$H \, \mathcal{I}[f] = 0 \mbox{.}$$
 \item[4.] Theorem \ref{T:main} includes a nonexistence result: there are no surfaces $f\left( \Sigma \right)$ satisfying \eqref{E:theeqn} and \eqref{E:smallness} with boundary conditions \ref{E:BCs} whose boundaries do not lie within a plane in $\mathbb{R}^3$.
\end{enumerate}

The structure of this article is as follows.  In Section \ref{S:prelim} we set up notation and state some fundamental results that are needed in the proof of Theorem \ref{T:main}.  In Section \ref{S:variation} we compute the normal variation of \eqref{E:energy} showing how \eqref{E:theeqn} and the boundary conditions \ref{E:BCs} arise.  In Section \ref{S:proof} we establish various estimates that culminate in the proof of Theorem \ref{T:main}.

\section{Preliminaries} \label{S:prelim}
\newtheorem{MSSobolev}{Theorem}[section]

Throughout this work we will employ cut-off functions defined as follows.  We take $\tilde{\gamma} \in C^2_c(\R^3)$ of the form $\tilde{\gamma}\left( x\right) = \hat{\gamma}\left( \frac{1}{\rho} \left| x\right| \right)$, for any $\rho>0$, where $\hat{\gamma} : \mathbb{R}_+\cup \left\{ 0 \right\} \rightarrow \left[ 0, 1\right]$ satisfies
$$\hat{\gamma}\left( s\right) = \begin{cases}
  1 & 0\leq s \leq \frac{1}{2}\\
  0 & s \geq 1 \mbox{.}
  \end{cases}$$
Then $\gamma=\tilde{\gamma}\circ f:\Sigma\rightarrow[0,1]$  satisfies
\begin{equation*}
\left\| \nabla\gamma \right\|_\infty \le c_{\gamma},\quad \mbox{ and }
\left\| \nabla^{2}\gamma \right\|_\infty \le c_{\gamma}(c_\gamma+|A|),
\end{equation*}
for $c_{\gamma} = \frac{c}{\rho}$ where $c>0$ is an absolute constant.

We will also need the Michael-Simon Sobolev inequality \cite{MS} for surfaces with boundary.  A proof in this setting appears for example in \cite{W14}.

\begin{MSSobolev}
	For $f:M^m \rightarrow \mathbb{R}^n$ a smooth immersion of $M$ with boundary $\partial M$ into $\mathbb{R}^n$ and any $u\in C^{1}\left( \overline{M}\right)$,
\begin{equation} \label{E:MS}
  \left[ \int_M \left| u\right|^{\frac{m}{m-1}} d\mu \right]^{\frac{m-1}{2}} \leq \frac{4^{m+1}}{\omega_m^{1/m}} \left[ \int_M \left( \left| \nabla u \right| + \left| H\right| \left| u \right| \right) d\mu + \int_{\partial M} \left| u\right| d\sigma \right]
  \end{equation}
	where $\omega_m$ is the volume of the unit ball in $\mathbb{R}^m$ and $d\sigma$ is the area element on $\partial M$.
\end{MSSobolev}

\noindent \textbf{Remark:} We only need to apply the above with $m=2$ and $n=3$ (or $n=k$ in view of the earlier Remark 2).  Moreover our boundary conditions \eqref{E:BCs} always ensure in our applications of the above that the boundary term  is identically equal to zero.  With these settings \eqref{E:MS} gives
$$\int_\Sigma u^2 d\mu \leq c \left[ \int_\Sigma \left( \left| \nabla u \right| + \left| H \right| \left| u \right| \right) d\mu \right]^2 \mbox{,}$$
where $c=\frac{32\sqrt{3}}{\sqrt{\pi}}$.

Let us finally mention that from the Codazzi equations 
$$\nabla_i h_{jk} = \nabla_j h_{ki} = \nabla_k h_{ij}$$
one can show (see, eg \cite{P17}) that the $k$-derivatives of the full $A$ tensor are controlled by those of $A^0$:
\begin{equation} \label{E:DkA}
\left| \nabla^{(k)} H \right| \leq \left| \nabla^{(k)} A \right| \leq c\left( n\right) \left| \nabla^{(k)} A^0 \right| \mbox{.}
\end{equation}

\section{Extremal surfaces with boundary for energy \eqref{E:energy}} \label{S:variation}
\newtheorem{var}{Lemma}[section]

We calculate the normal variation of energy \ref{E:energy} as follows.  In addition to previously-introduced notation we denote by $g_{ij}$ components of the metric on $f\left( \Sigma \right)$ and by $g^{ij}$ components of its inverse.

\begin{var}
	Given a smooth normal variation $\phi: \Sigma \rightarrow \mathbb{R}^3$ of $f: \Sigma \rightarrow \mathbb{R}^3$,
	\begin{multline} \label{E:variation}
	\left. \frac{d}{d\varepsilon}\mathcal{F}\left[ f+\epsilon\phi \right] \right|_{\epsilon=0} 
	=-2\int_{\Sigma} \left[ \Delta^{2}H+|A|^{2}\Delta H-(A^{0})^{ij}\nabla_{i}H\nabla_{j}H\right] \left< \phi, \nu \right> d\mu\\ 
	+2 \int_{\partial \Sigma} \left< \left( \Delta \phi + \left| A \right|^2 \phi\right) \nabla H + \nabla \Delta H \, \phi - \Delta H \, \nabla \phi, \eta\right> d\sigma \mbox{.}
	\end{multline}
\end{var}

Here $\nu$ denotes a smooth choice of unit normal to $f\left( \Sigma \right)$.

\noindent \textbf{Proof:} Writing $\varphi = \left< \phi, \nu \right>$, equation \eqref{E:variation} follows from the variations
$$\left. \frac{\partial}{\partial \epsilon}g^{\epsilon}_{ij}\right|_{\epsilon=0}=-2\varphi A_{ij},\;\; \left. \frac{\partial}{\partial \epsilon}g^{ij}_{\epsilon}\right|_{\epsilon=0}=2\varphi A^{ij},$$
$$\left. \frac{\partial}{\partial \epsilon}\sqrt{\det (g^{\epsilon}_{ij})}\right|_{\epsilon=0}=-H\varphi\sqrt{\det (g_{ij})}, \; \; \left. \frac{\partial}{\partial \epsilon}H_{\epsilon}\right|_{\epsilon=0}
=\Delta\varphi+\varphi|A|^{2}.$$
Calculations of these may be found in \cite{E} and \cite{MW}, for example.

Using the above we calculate, with slight abuse of notation and suppressing $\varepsilon$ where there is no chance of confusion
\begin{multline} \label{E:var1}
\left. \frac{d}{d\varepsilon} \int_\Sigma \left| \nabla H \right|^2 d\mu \right|_{\varepsilon=0}
=\left. \frac{d}{d\varepsilon} \int_U g^{ij} \frac{\partial}{\partial x_i} H \frac{\partial}{\partial x_j} H \sqrt{\det (g^{\epsilon}_{ij})} dx \right|_{\varepsilon=0}\\
=\int_\Sigma \left( 2 \varphi A^{ij} \right) \nabla_i H \nabla_j H d\mu + 2 \int_\Sigma \nabla^i \left( \Delta \varphi + \left| A\right|^2 \varphi \right) \nabla_j H d\mu + \int_\Sigma \left| \nabla H \right|^2 \left( -H \varphi \right) d\mu
 \end{multline}
 The first and last terms on the right hand side of \eqref{E:var1} combine to give the third term in \ref{E:variation}.  On the second term of \eqref{E:var1} we `integration by parts' on $\Sigma$ with boundary, that is, we apply the Divergence Theorem
$$\int_{\Sigma} \mbox{div}_\Sigma X \, d\mu = \int_{\partial \Sigma} \left< X, \eta \right> d\sigma \mbox{,}$$
where $X$ is a any smooth tangent vector field to $f\left( \Sigma \right)$, $\mbox{div}_\Sigma$ is the divergence on $f\left( \Sigma \right)$ and $\eta$ is the outer unit co-normal to $\partial \Sigma$.  (Of course a version also exists for $X$ a general vector field not necessarily tangent, for this we refer the reader to \cite{Si} for example.)

We have using the Divergence Theorem
$$\int_\Sigma \mbox{div} \left[ \left( \Delta\varphi+\varphi|A|^{2} \right) \nabla H \right] d\mu = \int_{\partial \Sigma} \left( \Delta\varphi+\varphi|A|^{2} \right) \left< \nabla H, \eta \right> d\sigma \mbox{;}$$
expanding out the left hand side by the product rule therefore yields
$$\int_\Sigma \nabla^i \left( \Delta\varphi+\varphi|A|^{2} \right) \nabla_i H  d\mu = - \int_\Sigma \left( \Delta\varphi+\varphi|A|^{2} \right) \Delta H d\mu + \int_{\partial \Sigma} \left( \Delta\varphi+\varphi|A|^{2} \right) \left< \nabla H, \eta \right> d\sigma  \mbox{.}$$
The second term of \eqref{E:variation} and the first boundary term are now clearly visible.  We integrate by parts twice more on the first term above:
\begin{multline*}
\int_\Sigma \Delta H \Delta \varphi \, d\mu = - \int_\Sigma \nabla^i \Delta H \nabla_i \varphi \, d\mu + \int_{\partial \Sigma} \Delta H \left< \nabla \varphi, \eta \right> d\mu\\
= \int_\Sigma \Delta^2 H \varphi \, d\mu -  \int_{\partial \Sigma} \left< \nabla \Delta H, \eta \right> \varphi \, d\mu 
+ \int_{\partial \Sigma} \Delta H \left< \nabla \varphi, \eta \right> d\mu 
\end{multline*}
revealing the remaining terms in \eqref{E:variation}.\hspace*{\fill}$\Box$
\mbox{}\\[8pt]

If $f\left( \Sigma \right)$ were closed without boundary, there would be no boundary terms in the expression \eqref{E:variation} and critical points of $\mathcal{F}\left[ f\right]$ would satisfy \eqref{E:theeqn}.

In our setting we impose \emph{flat boundary conditions} \eqref{E:BCs} on $\partial \Sigma$. Then the boundary terms in \eqref{E:variation} disappear (on the last term we integrate by parts on $\partial \Sigma$ that itself has no boundary) and we are left with \eqref{E:theeqn} for critical points of the energy.  We study smooth solutions \eqref{E:theeqn} with boundary conditions \eqref{E:BCs} and smallness condition \eqref{E:smallness}.

\section{Proof of Theorem \ref{T:main}} \label{S:proof}

We establish in turn estimates that facilitate the proof of Theorem \ref{T:main}.  The first four Lemmas below do not require the smallness condition \ref{E:smallness}.

Using the Divergence Theorem on $\Sigma$ (integration by parts) we begin with
\begin{lemma} \label{T:1}
	Surfaces satisfying \eqref{E:BCs} also satisfy
	\begin{align*}
	\int_\Sigma \left( \Delta H \right)^2 \gamma^p d\mu
	&= \int_\Sigma \mathcal{I}\left[ f\right] H \gamma^p d\mu + \int_\Sigma \left| A\right|^2 \left| \nabla H \right|^2 \gamma^p d\mu\\
	 & \quad  + \int_\Sigma H \nabla^i H \nabla_i \left| A \right|^2 \gamma^p d\mu 
	+ \int_\Sigma H \left( A^0\right)^{ij} \nabla_i H \nabla_j H \, \gamma^p d\mu\\
& \quad + p\int_\Sigma \left[ H\, \nabla^i \Delta H + \left( H \left| A \right|^2 - \Delta H \right) \nabla^i H \right] \nabla_i \gamma \cdot  \gamma^{p-1} d\mu
	\end{align*}
\end{lemma}

\noindent \textbf{Proof:} We multiply \eqref{E:theeqn} by $H \gamma^p$ for some constant $p>0$ to be chosen later:
\begin{multline} \label{E:zero}
		\int_\Sigma \mathcal{I}[f]\, H \gamma^{p} d\mu 
		= \int_\Sigma H \Delta^{2} H  \gamma^{p} d\mu + \int_\Sigma |A|^{2} H \Delta H \gamma^{p} d\mu\\ 
		- \int_\Sigma H \,\left( A^0\right)^{ij} \nabla_{i}H \nabla_{j} H \gamma^{p} d\mu \mbox{.}
		\end{multline}
Noting that 		
\begin{multline*}
	\int_\Sigma \mbox{div} \left( H \nabla\Delta H \gamma^{p}\right) d\mu 
	= \int_\Sigma H \Delta^{2} H \gamma^{p} d\mu + \int_\Sigma \nabla^i H \nabla_i\Delta H \gamma^{p} d\mu\\
	+ p\int_\Sigma H \nabla^i\Delta H \nabla_i \gamma \cdot \gamma^{p - 1} d\mu
= \int_{\partial \Sigma} H \nabla_{\eta} \Delta H   \gamma^{p}d\mu
= 0 
\end{multline*}
where the boundary integral is equal to zero in view of \eqref{E:BCs}, we have
\begin{equation} \label{E:one}
	\int_\Sigma H \Delta^{2} H\, \gamma^{p} d\mu = -\int_\Sigma \nabla^i H \nabla_i \Delta H \gamma^{p} d\mu - p\int_\Sigma H \nabla^i \Delta H \nabla_i\gamma \cdot \gamma^{p - 1} d\mu \mbox{.}
\end{equation}
Next we compute 
\begin{multline*}
\int_\Sigma \mbox{div} \left(\nabla H \Delta H \gamma^{p}\right) d\mu 
= \int_\Sigma \left( \Delta H\right)^{2} \gamma^{p} d\mu + \int_\Sigma \nabla^i H \nabla_i\Delta H \gamma^{p} d\mu \\+ p\int_\Sigma \gamma^{p - 1}\Delta H \nabla^i H \nabla_i \gamma \, d\mu
= \int_{\partial \Sigma} \Delta H \nabla_{\eta}H  \gamma^{p} d\mu
 = 0 \mbox{,}
\end{multline*}
where similarly \eqref{E:BCs} implies the boundary integral is equal to zero.  Hence
\begin{equation} \label{E:two}
-\int_\Sigma \nabla^i H \nabla_i\Delta H \gamma^{p} d\mu = \int_\Sigma \left( \Delta H\right)^{2} \gamma^{p} d\mu 
+ p\int_\Sigma \Delta H \nabla^i H \nabla_i\gamma \cdot \gamma^{p - 1} d\mu.
\end{equation}
Combining \eqref{E:one} and \eqref{E:two} we obtain
\begin{multline} \label{E:11}
\int_\Sigma H \Delta^{2} H \gamma^{p} d\mu\\
 = \int_\Sigma \left( \Delta H \right)^{2} \gamma^{p} d\mu + p\int_\Sigma \Delta H \nabla^i H \nabla_i \gamma\cdot \gamma^{p - 1} d\mu
- p\int_\Sigma H \nabla^i \Delta H \nabla_i\gamma\cdot \gamma^{p - 1} d\mu \mbox{.}
\end{multline}
Now
\begin{multline*}
	\int_\Sigma\mbox{div} \left(H|A|^{2}\nabla H\gamma^{p}\right)d\mu 
	= \int_\Sigma \nabla^i\left(|A|^{2}H\right) \nabla_i H \gamma^{p} d\mu + \int_\Sigma |A|^{2}H \Delta H \gamma^{p} d\mu \\
	+ p\int_\Sigma\gamma^{p - 1}H|A|^{2} \nabla^i H \nabla_i\gamma d\mu 
	= \int_{\partial \Sigma} H|A|^{2}\nabla_{\eta} H\gamma^{p} d\mu
	 = 0 
\end{multline*}
where we have again used \eqref{E:BCs}.  Thus
\begin{equation} \label{E:22}
\int_\Sigma |A|^{2}H \Delta H \gamma^{p} d\mu = - \int_\Sigma \nabla^i\left(|A|^{2}H\right) \nabla_i H \gamma^{p} d\mu - p\int_\Sigma\gamma^{p - 1} H |A|^{2}  \nabla^i H \nabla_i\gamma d\mu \mbox{.}
\end{equation}
Substituting \eqref{E:11} and \eqref{E:22} into \eqref{E:zero} we obtain
\begin{align*}
&\int_\Sigma \mathcal{I}[f] H \gamma^{p}d\mu \\
&= \int_\Sigma \left( \Delta H\right)^{2} \gamma^{p} d\mu + p \int_\Sigma\Delta H \nabla^i H \nabla_i \gamma \gamma^{p - 1} d\mu - p\int_\Sigma H \nabla^i\Delta H, \nabla_i\gamma \gamma^{p - 1} d\mu \\
& \quad - \int_\Sigma \nabla^i \left(|A|^{2}H\right) \nabla_i H \gamma^{p} d\mu - p\int_\Sigma \gamma^{p - 1}H|A|^{2} \nabla^i H \nabla_i{\gamma} d\mu 
- \int_\Sigma H  \left( A^0 \right)^{ij}\nabla_{i}H \nabla_{j} H \gamma^{p} d\mu \\
&\quad + p \int_\Sigma\left[\Delta H \nabla^i H \nabla_i\gamma - H \nabla^i\Delta H \nabla_i\gamma - H|A|^{2} \nabla^i H \nabla_i\gamma \right]\gamma^{p - 1} d\mu \mbox{.}
\end{align*}
The result follows.\hspace*{\fill}$\Box$

\begin{lemma} \label{T:2}
	Surfaces satisfying \eqref{E:BCs} also satisfy
	\begin{multline*}
	\int_\Sigma \left| \nabla^2 H \right|^2 \gamma^p d\mu
	 \leq c \int_\Sigma \mathcal{I}\left[ f\right] H \gamma^p d\mu + c \int_\Sigma \left| A\right|^2 \left| \nabla A \right|^2 \gamma^p d\mu + c \, c_\gamma^2 \int_\Sigma \left| \nabla A^0\right|^2 \gamma^{p-2} d\mu\\
+ p\int_\Sigma \left[ H\nabla^i \Delta H + \left( H \left| A \right|^2 - \Delta H \right) \nabla^i H \right] \nabla_i \gamma \cdot  \gamma^{p-1} d\mu
	\end{multline*}
\end{lemma}
\noindent \textbf{Proof:} Interchanging covariant derivatives, using the Codazzi equations and the Divergence theorem we may obtain exactly as in \cite{W14} that there is a universal constant $c$ such that 
\begin{multline*}
\frac{1}{c} \int_\Sigma \left(  \left| \nabla^2 H \right|^2 + H^2 \left| \nabla H \right|^2 \right) \gamma^p d\mu \\ \leq \int_\Sigma \left[ \left( \Delta H \right)^2  + \left| A^0\right|^2 \left| \nabla H \right|^2 \right] \gamma^p d\mu 
+ c_\gamma^2 \int_\Sigma \left| \nabla A^0 \right|^2 \gamma^{p-2} d\mu \mbox{.}
\end{multline*}
Further, we estimate 
\begin{equation*}
\int_\Sigma H \left( A^0\right)^{ij} \nabla_i H \nabla_j H \gamma^p d\mu
\leq \frac{1}{2c} \int_\Sigma H^2 \left| \nabla H \right|^2 \gamma^p d\mu + \frac{c}{2} \int_\Sigma \left| A^0 \right|^2 \left| \nabla H \right|^2 \gamma^p d\mu
\end{equation*}
and
$$\int_\Sigma H \, \nabla_i \left| A \right|^2 \nabla^i H \, \gamma^p d\mu \leq \tilde c \int \left| A\right|^2 \left| \nabla A \right|^2 \gamma^p d\mu \mbox{.}$$
Combining these with the Lemma \ref{T:1} yields the result.\hspace*{\fill}$\Box$

\begin{lemma} \label{T:3}
	Surfaces satisfying \eqref{E:BCs} also satisfy
	\begin{align*}
	&\int_\Sigma \left(  \left| \nabla^2 H \right|^2 + \left| A \right|^4 \left| A^0 \right|^2 \right) \gamma^p d\mu\\
	& \leq c \int_\Sigma \mathcal{I}\left[ f\right] H \gamma^p d\mu + c \int_\Sigma \left| A\right|^2 \left| \nabla A \right|^2 \gamma^p d\mu + c \int_\Sigma \left| A^0 \right|^6 \gamma^p d\mu + c \, c_\gamma^2 \int_\Sigma \left| \nabla A^0\right|^2 \gamma^{p-2} d\mu\\
	& \quad  + c\, c_\gamma^4 \int_\Sigma \left| A^0 \right|^2 \gamma^{p-4} d\mu + p\int_\Sigma \left[ H\nabla^i \Delta H + \left( H \left| A \right|^2 - \Delta H \right) \nabla^i H \right] \nabla_i \gamma \cdot  \gamma^{p-1} d\mu
	\end{align*}
\end{lemma}

\noindent \textbf{Proof:} By the same argument as for \cite[Lemma 5]{W14} we have
\begin{multline*}
\int_\Sigma \left( H^4 \left| A^0\right|^2 + H^2 \left| \nabla A^0\right|^2\right) \gamma^p d\mu\\
\leq c \int_\Sigma \left( H^2 \left| \nabla H\right|^2 +  \left|  A^0\right|^2 \left| \nabla A^0\right|^2 + \left| A^0\right|^6 \right) \gamma^p
+ c\, c_\gamma^4 \int_\Sigma \left| A^0 \right|^2 \gamma^{p-4} d\mu \mbox{.}
\end{multline*}
Moreover for $\varepsilon>0$,
\begin{multline*}
\int_\Sigma \left| A\right|^4 \left| A^0 \right|^2 \gamma^p d\mu = \int_\Sigma \left( H^2 + \frac{1}{2} \left| A^0 \right|^2 \right)^2 \left| A^0 \right|^2 d\mu\\
\leq \left( 1 + \varepsilon \right) \int_\Sigma H^4 \left| A^0 \right|^2 \gamma^p d\mu + \left( 1 + c\left( \varepsilon\right) \right) \int_\Sigma \left| A^0\right|^6 \gamma^p d\mu \mbox{.}
\end{multline*}
Thus from Lemma \ref{T:2} we obtain
\begin{align*}
&\int_\Sigma \left(  \left| \nabla^2 H \right|^2 + \left| A \right|^4 \left| A^0 \right|^2 \right) \gamma^p d\mu\\
&\leq c \int_\Sigma \mathcal{I}\left[ f\right] H \gamma^p d\mu + c \int_\Sigma \left| A\right|^2 \left| \nabla A \right|^2 \gamma^p d\mu + c \, c_\gamma^2 \int_\Sigma \left| \nabla A^0\right|^2 \gamma^{p-2} d\mu\\
&\quad + p\int_\Sigma \left[ H\nabla^i \Delta H + \left( H \left| A \right|^2 - \Delta H \right) \nabla^i H \right] \nabla_i \gamma \cdot  \gamma^{p-1} d\mu\\
&\quad + \left( 1 + \varepsilon \right) \int_\Sigma H^4 \left| A^0\right|^2 \gamma^p d\mu + \left( 1 + c\left( \varepsilon \right) \right) \int_\Sigma \left| A^0\right|^6 \gamma^p d\mu\\
&\leq c \int_\Sigma \mathcal{I}\left[ f\right] H \gamma^p d\mu + c \int_\Sigma \left| A\right|^2 \left| \nabla A \right|^2 \gamma^p d\mu\\
&\quad + \left( 1 + \varepsilon \right) \left[ c \int_\Sigma \left( H^2 \left| \nabla H\right|^2 +  \left|  A^0\right|^2 \left| \nabla A^0\right|^2 + \left| A^0\right|^6 \right) \gamma^p
+ c\, c_\gamma^4 \int_\Sigma \left| A^0 \right|^2 \gamma^{p-4} d\mu \right]\\
& \quad + c\, c_\gamma^2 \int_\Sigma \left| \nabla A^0 \right|^2 \gamma^{p-2} d\mu+ p\int_\Sigma \left[ H\nabla^i \Delta H + \left( H \left| A \right|^2 - \Delta H \right) \nabla^i H \right] \nabla_i \gamma \cdot  \gamma^{p-1} d\mu 
\end{align*}
which gives the statement of the Lemma. \hspace*{\fill} $\Box$

\begin{lemma} \label{T:4}
	Surfaces satisfying \eqref{E:BCs} also satisfy
	\begin{align*}
	&\int_\Sigma \left(  \left| \nabla^2 A \right|^2 + \left| A \right|^2 \left| \nabla A \right|^2 + \left| A \right|^4 \left| A^0 \right|^2 \right) \gamma^p d\mu\\
	&\leq c \int_\Sigma \mathcal{I}\left[ f\right] H \gamma^p d\mu + c \int_\Sigma \left| A\right|^2 \left| \nabla A \right|^2 \gamma^p d\mu + c \int_\Sigma \left| A^0 \right|^6 \gamma^p d\mu + c \, c_\gamma^2 \int_\Sigma \left| \nabla A^0\right|^2 \gamma^{p-2} d\mu\\
	& \quad  + c\, c_\gamma^4 \int_\Sigma \left| A^0 \right|^2 \gamma^{p-4} d\mu + p\int_\Sigma \left[ H\nabla^i \Delta H + \left( H \left| A \right|^2 - \Delta H \right) \nabla^i H \right] \nabla_i \gamma \cdot  \gamma^{p-1} d\mu
	\end{align*}
\end{lemma}

\noindent \textbf{Proof:} Using Simons' identity \cite{S68}
$$\Delta h_{ij} = \nabla_i \nabla_j H + H h_{im} h^{m}_{\ j} - \left| A\right|^2 h_{ij}$$
we may write
\begin{equation*}
  \Delta A^0_{ij}= \Delta h_{ij} - \frac{1}{2} g_{ij} \Delta H = \nabla_i \nabla_j H - \frac{1}{2} g_{ij} \Delta H + H h_{im} h^{m}_{\ j} - \left| A\right|^2 h_{ij} \mbox{.}
  \end{equation*}
  Using $h_{ij} = A^0_{ij} + \frac{1}{2} H g_{ij}$ this becomes
  \begin{multline*}
   \Delta A^0_{ij}=\nabla_i \nabla_j H - \frac{1}{2} g_{ij} \Delta H + H\left( A^0_{im} + \frac{1}{2} H g_{im}\right) \left( \left(A^0 \right)^m_{j} + \frac{1}{2} H g^m_{\ j} \right)\\
   = \nabla_i \nabla_j H - \frac{1}{2} g_{ij} \Delta H + H A^0_{im} \left( A^0\right)^m_{\ j} + \frac{1}{2} H^2 A^0_{ij} - \left| A^0\right|^2 A^0_{ij} - \frac{1}{2} H \left| A^0\right|^2 g_{ij} \mbox{.}
   \end{multline*}
   Hence for an absolute constant $c$ we have
   $$\left| \Delta A^0 \right| \leq \left| \nabla^2 H \right| + c \left| H \right| \left| A^0\right|^2 + c H^2 \left| A^0\right| + c \left| A^0\right|^3$$
   and so
   $$\left| \Delta A^0 \right|^2 \leq \left| \nabla^2 H \right|^2 + c H^4 \left| A^0\right|^2 + c\left| A^0 \right|^6 \mbox{.}$$
Interchange of second covariant derivatives and the Divergence Theorem then shows
\begin{multline*}
 \int_\Sigma \left| \nabla^2 A^0 \right|^2 \gamma^p d\mu\\
   \leq 2 \int_\Sigma \left| \Delta A^0 \right|^2 \gamma^p d\mu + c\int_\Sigma \left( \left| A\right|^2 \left| \nabla A\right|^2 + \left| A^0\right|^6 \right) \gamma^p d\mu
  + c\, c_\gamma^2 \int_\Sigma \left| \nabla A^0 \right|^2 \gamma^{p-2} d\mu \mbox{.}
\end{multline*}
Bearing in mind \eqref{E:DkA} and using also \cite{W14} inequality (31) we estimate
\begin{align*}
&\int_\Sigma \left(|\nabla^{2}A|^{2} + |A|^{2}|\nabla A|^{2} + |A|^{4}|A^{0}|^{2}\right) \gamma^{p} d\mu \\
&\leq 2\int_\Sigma \left| \Delta A\right|^{2}\gamma^{p} d\mu + c\int_\Sigma |A|^{2}|\nabla A|^{2} \gamma^{p} d\mu + \int_\Sigma|A|^{4}|A^{0}|^{2} \gamma^{p} d\mu + c\, c_{\gamma}^{2}\int_\Sigma|\nabla A^{0}|^{2}\gamma^{p - 2} d\mu \\
&\leq c\int_\Sigma |\nabla^{2}H|^{2} \gamma^{p} d\mu + c\int_\Sigma H^{4}|A^{0}|^{2}\gamma^{p}d\mu + \int_\Sigma|A^{0}|^{6}\gamma^{p}d\mu + c\int_\Sigma|A|^{2}|\nabla A|^{2} \gamma^{p} d\mu \\
& \quad + \int_\Sigma|A|^{4}|A^{0}|^{2} \gamma^{p} d\mu + c\, c_{\gamma}^2\int_\Sigma |\nabla A^{0}|^{2}\gamma^{p - 2} d\mu \\
\end{align*}
The result then follows using Lemma \ref{T:3}.\hspace*{\fill}$\Box$ 

\begin{lemma} \label{T:5}
	Surfaces satisfying \eqref{E:BCs} and \eqref{E:smallness} also satisfy
	\begin{equation*}
	\int_\Sigma \left(  \left| \nabla^2 A \right|^2 + \left| A \right|^2 \left| \nabla A \right|^2 + \left| A \right|^4 \left| A^0 \right|^2 \right) \gamma^p d\mu
	\leq c \int_\Sigma \mathcal{I}\left[ f\right] H \gamma^p d\mu 
	+ c\, c_\gamma^4 \int_\Sigma \left| A \right|^2 \gamma^{p-4} d\mu
	\end{equation*}
\end{lemma}

\noindent \textbf{Proof:}  Write $\left\| A \right\|_{2, [\gamma>0]}^2 = \int_{[\gamma>0]} \left| A \right|^2 d\mu$.  The idea is to use the smallness condition \eqref{E:smallness} to estimate the terms on the right hand side of Lemma \ref{T:4}.  In \cite{W14} it was shown using the Michael-Simon Sobolev inequality
\begin{multline*}
 \int_\Sigma \left( \left| A^0\right|^2 \left| A\right|^4 + \left| A \right|^2 \left| \nabla A\right|^2 \right) \gamma^p d\mu\\
 \leq c\left\| A \right\|_{2, [\gamma>0]}^2 \int_\Sigma \left(  \left| \nabla^2 A^0 \right|^2 + \left| A \right|^2 \left| \nabla A^0 \right|^2 + \left| A \right|^4 \left| A^0 \right|^2 \right) \gamma^p d\mu + c\, c_\gamma^4 \left\| A \right\|_{2, [\gamma>0]}^4 \mbox{;}
\end{multline*}
this result applies in the case of our boundary conditions \eqref{E:BCs}.  Thus we can absorb the non-$c_\gamma$ terms on the right hand side of Lemma \ref{T:4}.

We estimate the $c_\gamma$ terms from Lemma \ref{T:4} as follows:
$$c\, c_\gamma^4 \int_\Sigma \left| A^0 \right|^2 \gamma^{p-4} d\mu \leq c\, c_\gamma^4 \int_\Sigma \left| A \right|^2 \gamma^{p-4} d\mu \mbox{;}$$
via the Divergence Theorem and the Cauchy-Schwarz and Peter-Paul inequalities we have for $\varepsilon>0$
$$c\, c_\gamma^2 \int_\Sigma \left| \nabla A^0 \right|^2 \gamma^{p-2} d\mu \leq \varepsilon \int_\Sigma \left| \nabla^2 A\right|^2 \gamma^p d\mu + c\left( \varepsilon\right) c_\gamma^4 \int_\Sigma \left| A \right|^2 \gamma^{p-4} d\mu \mbox{;}$$
with this in turn we estimate
\begin{multline*}
 p \int_\Sigma \Delta H \nabla^i H \nabla_i \gamma \cdot \gamma^{p-1} d\mu\\
  \leq c\,  c_\gamma \int_\Sigma \left| \nabla^2 A \right| \left| \nabla H \right| \gamma^{p-1} d\mu
  \leq \varepsilon \int_\Sigma \left| \nabla^2 A\right|^2 \gamma^p d\mu + c\left( \varepsilon\right) c_\gamma^4 \int_\Sigma \left| A \right|^2 \gamma^{p-4} d\mu
\end{multline*}
and
\begin{equation*}
p \int_\Sigma H \left| A \right|^2 \nabla^i H \nabla_i \gamma \cdot \gamma^{p-1} d\mu 
\leq \varepsilon \int_\Sigma \left| A \right|^2 \left| \nabla A^0 \right|^2 \gamma^p d\mu + c\left( \varepsilon\right) c_\gamma^2 \int_\Sigma H^2 \left| A \right|^2 \gamma^{p-2} d\mu \mbox{.}
\end{equation*}
Now by the Michael-Simon Sobolev inequality ($\left| A \right|=0$ on $\partial \Sigma$ from \eqref{E:BCs})
\begin{align*}
\int_\Sigma H^2 \left| A \right|^2 \gamma^{p-2} d\mu &\leq \int_\Sigma \left| A \right|^4 \gamma^{p-2} d\mu \\
&\leq c \left( \int_\Sigma \left| \nabla \left| A \right|^2 \right| \gamma^{\frac{p-2}{2}} d\mu \right)^2 +  \left( \int_\Sigma \left| A \right|^3 \gamma^{\frac{p-2}{2}} d\mu \right)^2\\
& \leq c  \left( \int_\Sigma \left| \nabla A \right| \left| A \right| \gamma^{\frac{p-2}{2}} d\mu \right)^2 + c\, c_\gamma^2 \left( \int_\Sigma \left| A \right|^2 \gamma^{\frac{p-4}{2}} d\mu \right)^2\\
& \quad  + c \left\| A \right\|_{2, [\gamma>0]}^2 \int_\Sigma \left| A\right|^4 \gamma^{p-2} d\mu \mbox{.}
\end{align*}
Absorbing on the left and using the Cauchy-Schwarz inequality we obtain
\begin{equation*}
\int_\Sigma  \left| A \right|^4 \gamma^{p-2} d\mu
\leq c \left\| A \right\|_{2, [\gamma>0]}^2 \int_\Sigma \left| \nabla A\right|^2 \gamma^{p-2} d\mu + c\, c_\gamma^2 \left\| A \right\|_{2, [\gamma>0]}^4
\end{equation*}
and so
\begin{equation*}
c\, c_\gamma^2 \int_\Sigma H^2 \left| A \right|^2 \gamma^{p-2} d\mu\\
\leq c \, c_\gamma^2 \left\| A \right\|_{2, [\gamma>0]}^2 \int_\Sigma \left| \nabla A\right|^2 \gamma^{p-2} d\mu + c\, c_\gamma^4 \left\| A \right\|_{2, [\gamma>0]}^4 \mbox{.}
\end{equation*}

%

%
%
For the remaining term from Lemma \ref{T:4}  we use the Divergence Theorem ($H=0$ on $\partial \Sigma$ in view of \eqref{E:BCs})
\begin{multline*}
\int_\Sigma H \nabla^i \Delta H \nabla_i \gamma \cdot \gamma^{p-1} d\mu 
= - \int_\Sigma \Delta H \nabla^i H \nabla_i \gamma \cdot \gamma^{p-1} d\mu - \int_\Sigma H \Delta H \Delta \gamma \cdot \gamma^{p-1} d\mu\\
 - \left( p-1\right) \int_\Sigma H \Delta H \left| \nabla \gamma \right|^2 \gamma^{p-2} d\mu \mbox{.}
\end{multline*}
We now estimate for $\varepsilon>0$,
\begin{equation*}
- \int_\Sigma \Delta H \nabla^i H \nabla_i \gamma \cdot \gamma^{p-1} d\mu
\leq \varepsilon \int_\Sigma \left( \Delta H \right)^2 \gamma^p d\mu + c\left( \varepsilon\right) c_\gamma^2 \int_\Sigma \left| \nabla H \right|^2 \gamma^{p-2} d\mu \mbox{;}
\end{equation*}
\begin{align*}
&- \int_\Sigma H \Delta H \Delta \gamma \cdot \gamma^{p-1} d\mu\\
& \leq c\, c_\gamma \int_\Sigma \left| H\right| \left| \Delta H \right| \left( c_\gamma + \left| A \right| \right) \gamma^{p-1} d\mu\\
& \leq \varepsilon \int_\Sigma \left( \Delta H \right)^2 \gamma^p d\mu + c\, c_\gamma^4 \int_\Sigma \left| A\right|^2 \gamma^{p-2} d\mu + c\, c_\gamma^2 \int_\Sigma H^2 \left| A \right|^2 \gamma^{p-2} d\mu 
\end{align*}
and
\begin{equation*}
- \left( p-1\right) \int_\Sigma H \Delta H \left| \nabla \gamma \right|^2 \gamma^{p-2} d\mu
\leq \varepsilon  \int_\Sigma \left( \Delta H \right)^2 \gamma^p d\mu + c\, c_\gamma^4 \int_\Sigma H^2 \gamma^{p-4} d\mu
\end{equation*}
Inserting all these estimates and absorbing on the left yields the result.\hspace*{\fill}$\Box$
\mbox{}\\[8pt]

\noindent \textbf{Completion of the proof of Theorem \ref{T:main}:} Using Lemma \ref{T:5}, surfaces satisfying \eqref{E:BCs}, \eqref{E:smallness} and \eqref{E:theeqn} also satisfy
$$\int_\Sigma \left(  \left| \nabla^2 A \right|^2 + \left| A \right|^2 \left| \nabla A \right|^2 + \left| A \right|^4 \left| A^0 \right|^2 \right) \gamma^p d\mu
	\leq \frac{c}{\rho^4} \int_\Sigma \left| A \right|^2 \gamma^{p-4} d\mu \leq \frac{c}{\rho^4} \varepsilon_0$$
for an absolute constant $c$.  With $p=4$, taking $\rho\rightarrow \infty$ we see that $f\left( \Sigma\right)$ must have
$$\left| A\right|^4 \left| A^0 \right|^2 \equiv 0 \mbox{.}$$
Since 
$$\left| A^0 \right|^6 \leq \left| A\right|^4 \left| A^0\right|^2$$
 we have that
  $$\left| A^0 \right|^2 \equiv 0$$
  implying $f\left( \Sigma \right)$ is either part of a sphere or part of a plane.  The boundary condition \eqref{E:BCs} implies $f\left( \Sigma\right)$ is part of a plane.\hspace*{\fill}$\Box$


\end{document}